\documentclass[11pt,a4paper]{article}
\usepackage[T1]{fontenc}
\usepackage{lmodern}
\usepackage[utf8]{inputenc}
\usepackage[english]{babel}
\usepackage[margin=28mm]{geometry}
\usepackage{amsmath,amssymb,amsthm,mathtools}
\usepackage[hidelinks]{hyperref}
\usepackage[expansion=false]{microtype}

\newcommand{\End}{\operatorname{End}}
\newcommand{\G}{{\mathfrak{g}}}
\newcommand{\h}{{\mathfrak{h}}}
\newcommand{\ad}{\mathrm{ad}}

\newcommand{\esp}{\quad\mbox{and}\quad}
\newcommand{\R}{\mathbb{R}}
\newcommand{\md}{\mathrm d}

\newcommand{\Id}{\operatorname{Id}}

\newcommand{\im}{\mathrm{Im}}

\hypersetup{
	pdftitle={Left-Invariant Kundt Structures on Almost Abelian Lie Groups},
	pdfauthor={S. El Bourkadi and M. W. Mansouri}
}

\newtheorem{theorem}{Theorem}[section]
\newtheorem{pro}[theorem]{Proposition}
\newtheorem{Lem}[theorem]{Lemma}
\newtheorem{corollary}[theorem]{Corollary}
\theoremstyle{definition}
\newtheorem{Def}[theorem]{Definition}
\newtheorem{example}[theorem]{Example}
\theoremstyle{remark}
\newtheorem{remark}[theorem]{Remark}
\newcommand{\spn}{\operatorname{span}}

\numberwithin{equation}{section}
\allowdisplaybreaks
\title{Left-Invariant Kundt Structures on Almost Abelian Lie Groups}

\author{S. El Bourkadi and M. W. Mansouri\\Universit\'e Ibn Tofail\\Facult\'e des Sciences, Laboratoire L.A.G.A.\\K\'enitra, Morocco\\\texttt{mansourimohammed.wadia@uit.ac.ma}\\\texttt{said.elbourkadi@uit.ac.ma}}
\date{}

\begin{document}
	\maketitle
	
	\begin{abstract}
		We study left-invariant Kundt structures on almost abelian Lie groups.
		We establish algebraic existence criteria and characterize the admissible
		real Jordan forms according to the signature of the Abelian ideal.
		We then characterize the Kundt pairs defining pp-waves and plane waves
		in arbitrary dimension, with particular attention to recurrent and
		parallel null directions. Examples in dimension four illustrate both
		plane and non-plane pp-wave geometries.
	\end{abstract}
	
	\noindent\textbf{Keywords.} Kundt structures, almost abelian Lie groups,
	Lorentzian geometry, pp-waves, plane waves, null vector fields.

	\smallskip
	\noindent\textbf{2020 Mathematics Subject Classification.}
	Primary 53C50; Secondary 22E25, 53Z05, 17B30.
	\section{Introduction}

	The geometry of Lie groups endowed with left-invariant Lorentzian metrics provides a natural framework for the study of homogeneous spacetimes. Among the geometric structures of particular interest are Kundt structures, which are characterized by the existence of a distinguished null direction satisfying suitable geometric conditions. These structures arise naturally in Lorentzian geometry and mathematical physics and have attracted considerable attention in recent years.
	
	Boucetta, Meliani and Zeghib developed an algebraic approach to left-invariant Kundt structures and investigated their classification on three-dimensional Lie groups \cite{B-M-Z}. Their characterization relates the Kundt condition to the existence of a degenerate codimension-one Lie subalgebra preserved by the Levi-Civita product. Further developments concerning Kundt geometry, null distributions and their geometric properties can be found in \cite{M-Z,A}.
	
	The aim of the present paper is to investigate left-invariant Kundt structures on the class of almost abelian Lie groups. An almost abelian Lie algebra is a non-Abelian Lie algebra admitting an Abelian ideal of codimension one. Such Lie algebras have been extensively studied due to their relatively simple algebraic structure and the richness of the geometric structures they support \cite{Av}. In particular, every almost abelian Lie algebra can be written as
	\[
	\mathfrak g_A=V\rtimes_A\mathbb Re,
	\qquad [e,v]=Av,\quad [V,V]=0,
	\]
	where \(A\) is a nonzero real endomorphism of \(V\). This description allows geometric questions to be translated into problems involving linear algebra and the spectral properties of \(A\).
	
	The class of almost abelian Lie groups contains several important examples arising in geometry and mathematical physics. In low dimensions, these include the group of Euclidean motions \(E(2)\), the group of isometries of the Minkowski plane \(E(1,1)\), and the affine group \(\operatorname{Aff}(\mathbb R)\). These examples illustrate the diversity of algebraic and geometric properties within this class, particularly with respect to unimodularity, curvature and geodesic completeness.
	
	The study of left-invariant Lorentzian metrics on almost abelian Lie groups has recently attracted considerable attention. In particular, Bajo, Benayadi and Lebzioui \cite{B-B} investigated flat left-invariant Lorentzian metrics on almost abelian Lie groups, providing existence and classification results for such geometric structures. 
	
	More recently, Sato and Tsuyuki \cite{S-T} studied Ricci-flat left-invariant Lorentzian metrics on almost abelian Lie groups in arbitrary dimension, including a generalization of the Petrov solution. They established algebraic conditions for Ricci-flatness and identified certain families of plane waves. These investigations provide a useful background for the present work. However, the Kundt condition is more general and does not require flatness, Ricci-flatness, unimodularity or nilpotency. Our approach is therefore different: we focus on the algebraic characterization of left-invariant Kundt structures without imposing curvature restrictions.
	
	In this paper, we establish algebraic criteria for the existence of left-invariant Kundt structures on almost abelian Lie groups. We first show that every almost abelian Lie algebra admits a Kundt structure and distinguish two types of Kundt pairs according to the position of the associated codimension-one subalgebra. We then obtain necessary and sufficient conditions in terms of the defining endomorphism \(A\), taking into account the possible signatures of the metric restricted to the Abelian ideal. In particular, we derive explicit characterizations involving real Jordan normal forms.
	
	We further investigate recurrent and parallel null directions and characterize the Kundt structures that define pp-waves and plane waves in arbitrary dimension. This extends the study beyond the three-dimensional setting and complements recent work on four-dimensional pp-wave Lie groups and harmonic curvature \cite{G-R-V}. Explicit examples in dimension four illustrate the distinction between plane-wave and non-plane-wave geometries. The relationship between locally homogeneous pp-waves
	and plane waves has been investigated by Globke and
	Leistner \cite{G}.
	
	The paper is organized as follows. In Section 2, we recall the algebraic characterization of left-invariant Kundt structures and establish preliminary results on almost abelian Lie algebras. In Section 3, we investigate Kundt pairs of types (I) and (II), obtaining necessary and sufficient existence conditions and characterizations in terms of real Jordan normal forms. Finally, in Section 4, we characterize Kundt structures of pp-wave and plane-wave type in arbitrary dimension, providing explicit criteria and illustrative examples.

	\section{Algebraic Framework}
	In this section, we recall the algebraic characterization
	of left-invariant Kundt structures and establish the
	preliminary results needed for their study on almost
	abelian Lie groups.

	\subsection{Lorentzian Lie algebras and Kundt pairs}
	
	Let $\G$ be an $n$-dimensional real Lie algebra endowed with a Lorentzian inner product $\langle\cdot,\cdot\rangle$, that is, a non-degenerate symmetric bilinear form of signature $(1,n-1)$.

	The corresponding Levi-Civita product $
	\nabla:\G\times\G
	\longrightarrow\G$
	is uniquely determined by the Koszul formula
	\begin{equation}	\label{koszul}
		2\langle\nabla_x y,z\rangle
		=\langle[x,y],z\rangle
		-\langle[y,z],x\rangle
		+\langle[z,x],y\rangle,
	\end{equation}
	for all $x,y,z\in\G$.
	
	Following the algebraic framework introduced by  Boucetta, Meliani, and Zeghib \cite{B-M-Z} (see also \cite{M}), we recall the definition of a Kundt pair. For more details on Kundt structures, see \cite{M-Z}.
	
	\begin{Def}
		Let $(\G,\langle\cdot,\cdot\rangle)$ be a Lorentzian
		Lie algebra endowed with its Levi-Civita product $\nabla$.
		A \emph{Kundt pair} on $\G$ is a pair
		$(\langle\cdot,\cdot\rangle,\h)$, where
		$\mathfrak h$ is a degenerate Lie subalgebra of codimension-one	satisfying the following conditions:
		\begin{enumerate}
			\item $\h$ is $\nabla$-invariant, i.e.\ $\nabla_x y\in\mathfrak h$
			for all $x,y\in\h$.
			\item For any nonzero vector $p$ such that	$\h^\perp=\R p$, one has $\nabla_p p=0$.
		\end{enumerate}
	\end{Def}
	A pair $(\langle\cdot,\cdot\rangle,\mathfrak h)$
	satisfying the first condition is called a
	\emph{locally Kundt pair}.
	It is called a \emph{Kundt pair} if, in addition,
	$\nabla_p p=0$.	Recall that a nonzero vector $p\in\mathfrak g$ is said to be \emph{null (or lightlike)} if $\langle p,p\rangle=0$.
	Since $\mathfrak h$ is degenerate and of codimension-one,
	its orthogonal complement is a null line
	$\mathfrak h^\perp=\mathbb Rp\subseteq\mathfrak h$.
	The additional requirement $\nabla_p p=0$ ensures that
	the integral curves of the corresponding left-invariant
	null vector field are affinely parametrized null geodesics,
	thereby distinguishing Kundt pairs from locally Kundt pairs.

	However, there exists a large class of Lie groups that cannot carry any left-invariant locally Kundt structure. For example, compact semisimple Lie groups admit no left-invariant locally Kundt structure; see \cite{B-M-Z} for more details.

	\begin{Def}
		Let $(\G,\langle\cdot,\cdot\rangle)$ be a
		Lorentzian Lie algebra, and let
		$p\in\G$ be a nonzero isotropic vector.	The vector $p$ is called
		\emph{recurrent} if $\nabla_x p\in\R p,$
		$ \forall x\in\G$,
		and \emph{parallel} if $\nabla_x p=0$,
		$ \forall x\in\G$.
		
	\end{Def}
	
	Every parallel isotropic vector is recurrent, but the converse need not hold.

	The following Lemma avoids	computing all the coefficients of the	Levi-Civita connection.
	\begin{Lem}	\label{k-c}
		Let	$(\G,\langle\cdot,\cdot\rangle)$ be a Lorentzian Lie algebra, and let	$\h$ be a degenerate codimension-one Lie subalgebra.	Suppose $\h^\perp=\R p$, $p\neq0$.	Then $(\langle\cdot,\cdot\rangle,\h)$	is a Kundt pair if and only if
		\begin{equation}
			\langle[p,x],y\rangle+\langle x,[p,y]\rangle=0,
			\quad\forall x,y\in\h,
			\label{2.2}
		\end{equation}
		and
		\begin{equation}
			\langle[p,z],p\rangle=0,\quad\forall z\in\G.
			\label{2.3}
		\end{equation}
	\end{Lem}
	
	\begin{proof}
		Since $\h$ is a Lie subalgebra and	$\h=p^\perp$, the Koszul formula yields	\[2\langle\nabla_x y,p\rangle=\langle[p,x],y\rangle			+\langle x,[p,y]\rangle,\qquad\forall x,y\in\h.\]
		Thus, $\nabla_{\h}\h\subseteq\h$ is equivalent	to \eqref{2.2}.
		Moreover, for all $z\in\G$, $	2\langle\nabla_p p,z\rangle
		=-2\langle[p,z],p\rangle$.
		
		By nondegeneracy of the metric,
		$\nabla_pp=0$ is equivalent to
		\eqref{2.3}.
	\end{proof}
	
	\subsection{Almost abelian Lie algebras}
	
	In this section, we recall some basic algebraic properties  of real almost abelian	Lie algebras of finite dimension.
	\begin{Def}
		A non-abelian real Lie algebra $\G$ is called \emph{almost abelian} if it admits an	abelian ideal of codimension-one.
	\end{Def}
	Fix such an ideal $V$. Then $\G_A=V\rtimes_A\R e$,	where $	A=\operatorname{ad}_e|_V\in\End(V)$, $A\neq0$. The Lie bracket is given by
	\begin{equation}
		[v+te,w+se]=tAw-sAv,\quad \forall v,w\in V\quad\forall  s,t\in\R
	\end{equation}
	In particular, $[\G_A,\G_A]	=\im A$, and $Z(\G_A)=\ker A$.	The linear endomorphism $A$ therefore encodes the Lie algebra structure	of $\G_A$.	Throughout this paper, an almost abelian Lie algebra is assumed to be non-abelian, that is, \(A\neq 0\).
	
	The following proposition classifies the codimension-one subalgebras of an almost abelian Lie algebra and will play a crucial role in our study of Kundt structures.
	
	\begin{pro}	\label{Pr 2-5}
		Let $\G_A=V\rtimes_A\R e$ 	be an almost abelian Lie algebra.	Every codimension-one Lie subalgebra $\h$ of $\G_A$ is of	one of the following two types:
		\begin{enumerate}
			\item[\textup{(I)}] $\h=V$;
			\item[\textup{(II)}] $\h=W\oplus\mathbb R(e+v_0)$, where $W\subset V$  is  an $A$-invariant hyperplane of $V$  and  $v_0\in V$.
		\end{enumerate}
		Moreover, any subalgebra of type \textup{(II)} is mapped onto
		$\h=W\oplus\R e$ by some automorphism of $\G_A$.
	\end{pro}
	
	\begin{proof}
		If $\h\subseteq V$, then $\h=V$ by dimension.	Otherwise, choose $v_0\in V$ such that $e+v_0\in\h$ and set $W=\h\cap V$.	Then $W$ is a hyperplane of $V$ and $\h=W\oplus\R(e+v_0)$.	Since $\h$ is a Lie subalgebra, for every	$w\in W$,
		\[Aw=[e+v_0,w]\in\h\cap V=W.\]
		Thus $A(W)\subseteq W$. The converse is immediate.
		
		Finally, the map $\Phi(e)=e-v_0$, $\Phi|_V=\mathrm{Id}_V$ is an automorphism of $\G_A$ satisfying $\Phi(\h)=W\oplus\R e$.	This completes the proof.
	\end{proof}

\section{Kundt Structures on Almost Abelian Lie Algebras}\label{k-a-a}

Let $\G_A=V\rtimes_A\R e$ be an almost abelian Lie algebra endowed with
a Lorentzian inner product $\langle\cdot,\cdot\rangle$.
By Proposition~\ref{Pr 2-5}, a Kundt pair is, up to an automorphism of
$\G_A$, of one of the following two types:
\[
\textup{(I)}\quad \h=V,
\qquad
\textup{(II)}\quad \h=W\oplus\R e,\quad A(W)\subseteq W,
\]
where $W$ is a hyperplane of $V$. The automorphism of
Proposition~\ref{Pr 2-5} fixes $V$ pointwise, and the inner product is
transported by it; hence the signature of $\langle\cdot,\cdot\rangle|_V$
is not affected. This signature (Euclidean, degenerate or Lorentzian)
will serve to organize the results.

We say that $\h$ is a \emph{Kundt subalgebra} if
$(\langle\cdot,\cdot\rangle,\h)$ is a Kundt pair. In particular
$\h$ is degenerate and $\h^\perp=\R p\subseteq\h$.
\subsection{Kundt subalgebras of type (I)}
The following theorem characterizes Kundt structures of type (I).
\begin{theorem}\label{typeI}
	Let $\langle\cdot,\cdot\rangle$ be a Lorentzian inner product on $\G_A$.
	Then $(\langle\cdot,\cdot\rangle,V)$ is a Kundt pair if and only if $V$
	is degenerate. In this case $V^\perp=\R p\subset V$, and there are a
	Euclidean complement $E$ of $\R p$ in $V$ and a choice of $e$ and $p$ such that
	\begin{equation}\label{gram-I}
		\langle e,p\rangle=1,\qquad
		\langle e,e\rangle=\langle p,p\rangle=0,\qquad
		E\perp\spn\{e,p\},
	\end{equation}
	and
	\begin{equation}\label{L-B}
		[e,p]=Ap=ap+b,\qquad
		[e,x]=Ax=\alpha(x)p+Bx,\quad\forall x\in E,
	\end{equation}
	for some $a\in\R$, $b\in E$, $\alpha\in E^*$, $B\in\End(E)$ are not all zero.
	
	Conversely, let $(E,\langle\cdot,\cdot\rangle_E)$ be a Euclidean vector
	space, $V=\R p\oplus E$, $\G_A=\R e\oplus V$, and let $a\in\R$, $b\in E$,
	$\alpha\in E^*$, $B\in\End(E)$. Then \eqref{L-B} defines an almost
	abelian Lie algebra, \eqref{gram-I} defines a Lorentzian inner product on
	it, and $(\langle\cdot,\cdot\rangle,V)$ is a Kundt pair.
\end{theorem}

\begin{proof}
	A Kundt subalgebra is degenerate, which gives the ``only if'' part.
	Conversely, let $V$ be degenerate. Since $V$ is a hyperplane,
	$V^\perp=\R p$ is contained in $V$. As $V$ is abelian, $[p,x]=0$ for all
	$x\in V$, so \eqref{2.2} holds. As $V$ is an ideal, $[p,z]\in V$ for all
	$z\in\G_A$, and $p\perp V$ gives \eqref{2.3}. Lemma~\ref{k-c} applies.
	
	Normal form. The form induced on $V/\R p$ is positive definite, so any
	complement $E$ of $\R p$ in $V$ is Euclidean. Since the inner product is
	nondegenerate and $p\perp V$, we have $\langle e,p\rangle\neq0$, and we
	rescale $p$ so that $\langle e,p\rangle=1$. Replacing $e$ by $e+v$ with
	$v\in V$ does not change $A=\ad_e|_V$ ($V$ is abelian). Choose first
	$x\in E$ with $\langle e+x,y\rangle=0$ for all $y\in E$, then
	$\lambda\in\R$ such that $e+x+\lambda p$ is null. The new $e$ satisfies
	\eqref{gram-I}. Decomposing $Ap$ and $Ax$ along $V=\R p\oplus E$ gives
	\eqref{L-B}.
	
	Converse. The brackets \eqref{L-B} define an almost abelian Lie algebra
	with abelian ideal $V$. In the basis $(e,p,E)$ the Gram matrix is
	$\left(\begin{smallmatrix}0&1&0\\1&0&0\\0&0&I\end{smallmatrix}\right)$,
	so the inner product is Lorentzian and $V^\perp=\R p$; hence $V$ is
	degenerate, and the first part applies.
\end{proof}

\begin{corollary}\label{cor-exist-I}
	Every almost abelian Lie algebra admits a Kundt pair of type~\textup{(I)}.
\end{corollary}

\begin{proof}
	Choose $0\neq p\in V$ and a complement $E$ of $\R p$ in $V$, and endow $E$
	with a Euclidean inner product. Every $A\in\End(V)$ has the form
	\eqref{L-B} for suitable $a,b,\alpha,B$. The converse part of
	Theorem~\ref{typeI} concludes.
\end{proof}

The following criterion describes when the null line $\R p$ of a Kundt
pair of type~\textup{(I)} is recurrent or parallel; it will be used in
Section~\ref{sec:ppwaves}.

\begin{pro}\label{p-recurrent}
	Let $(\langle\cdot,\cdot\rangle,V)$ be a Kundt pair, $V^\perp=\R p$. Then
	$p$ is recurrent if and only if $Ap=\lambda p$ for some $\lambda\in\R$.
	In this case
	\[
	\nabla_e p=\lambda p,\qquad \nabla_v p=0\quad\forall v\in V.
	\]
	Moreover, $p$ is parallel if and only if $Ap=0$.
\end{pro}

\begin{proof}
	Since $V$ is abelian and an ideal, the Koszul formula gives, for
	$v,w\in V$,
	\[
	2\langle\nabla_v p,w\rangle=0,\qquad
	2\langle\nabla_v p,e\rangle=\langle Ap,v\rangle,\qquad
	2\langle\nabla_e p,v\rangle=\langle Ap,v\rangle,\qquad
	\langle\nabla_e p,e\rangle=\langle Ap,e\rangle,
	\]
	where we used $\langle Av,p\rangle=0$. In particular
	$\nabla_v p\in V^\perp=\R p$ for all $v\in V$.
	
	If $p$ is recurrent, then $\nabla_e p\in\R p\perp V$, hence
	$\langle Ap,v\rangle=0$ for all $v\in V$, i.e.\ $Ap\in V^\perp=\R p$.
	Conversely, let $Ap=\lambda p$. Then $\langle Ap,v\rangle=0$ for all
	$v\in V$, so $\nabla_v p\perp V\oplus\R e$ gives $\nabla_v p=0$, while
	$\nabla_e p\in V^\perp=\R p$ and $\langle\nabla_e p,e\rangle=\lambda\langle p,e\rangle$
	give $\nabla_e p=\lambda p$. Finally, $p$ is parallel if and only if it is
	recurrent with $\lambda=0$.
\end{proof}

\subsection{Kundt subalgebras of type (II)}\label{subsec:typeII}

Let $W\subset V$ be an $A$-invariant hyperplane and assume that
$\h=W\oplus\R e$ is degenerate, with $\h^\perp=\R p\subseteq\h$. Two cases
occur according to the position of the null generator $p$ with respect
to $W$.

\begin{theorem}\label{th-typeII}
	With the above notation:
	\begin{enumerate}
		\item[\textup{(II-a)}] If $p\notin W$, normalize $p=e+u$ with $u\in W$.
		Then $(\langle\cdot,\cdot\rangle,\h)$ is a Kundt pair if and only if
		\[
		\im A\subseteq W
		\esp
		\langle Aw_1,w_2\rangle+\langle w_1,Aw_2\rangle=0,\quad
		\forall\, w_1,w_2\in W.
		\]
		In this case $W$ is Euclidean.
		\item[\textup{(II-b)}] If $p\in W$, then $(\langle\cdot,\cdot\rangle,\h)$ is
		a Kundt pair if and only if $Ap\in\R p$. In this case $V$ is Lorentzian.
	\end{enumerate}
\end{theorem}

\begin{proof}

	(II-a) Let $v\in V$ we have $[p,v]=Av$, and $[p,e]=-Au$. Condition
	\eqref{2.2} on $W\times W$ says that $A|_W$ is skew-adjoint. The other
	instances of \eqref{2.2} then hold automatically. For $w\in W$, since
	$Aw\in W\subset p^\perp$ we have $\langle e,Aw\rangle=-\langle u,Aw\rangle$,
	so
	\[
	\langle[p,w],e\rangle+\langle w,[p,e]\rangle
	=-\bigl(\langle Aw,u\rangle+\langle w,Au\rangle\bigr)=0;
	\]
	and $\langle[p,e],e\rangle=-\langle Au,e\rangle=\langle Au,u\rangle=0$
	by skew-adjointness, since $Au\in W$. Next, \eqref{2.3} reads
	$\langle Av,p\rangle-t\langle Au,p\rangle=0$ for all $v+te$, that is
	$\langle Av,p\rangle=0$ for all $v\in V$; equivalently
	$\im A\subseteq V\cap p^\perp=V\cap\h=W$. Finally $\h$ is a degenerate
	hyperplane of a Lorentzian space and $W$ is a complement of its radical
	$\R p$, so $W\simeq\h/\R p$ is Euclidean.
	
	(II-b) Here $[p,W]=0$ and $[p,e]=-Ap$. Condition \eqref{2.2} is trivial on
	$W\times W$; it reads $\langle w,Ap\rangle=0$ for $w\in W$ and
	$\langle e,Ap\rangle=0$ otherwise. Hence \eqref{2.2} is equivalent to
	$Ap\perp\h$, i.e.\ $Ap\in\h^\perp=\R p$. Then \eqref{2.3} is obvious,
	because $[p,z]\in\R p$ for all $z$ and $p$ is null.
	It remains to see that $V$ is Lorentzian. It cannot be Euclidean, since
	it contains the null vector $p$. If it were degenerate, with
	$V^\perp=\R q$, then the null vector $p\in V$ would lie in $V^\perp$
	(the form induced on $V/V^\perp$ is positive definite), so $p\in\R q$;
	but $p\perp e$ whereas $\langle q,e\rangle\ne0$ by nondegeneracy, a
	contradiction.
\end{proof}

We now translate Theorem~\ref{th-typeII} into existence criteria on $A$.
Consider the following conditions.
\begin{enumerate}
	\item[\textup{(a)}] There is an $A$-invariant flag $L\subset W\subset V$ with
	$\dim L=1$ and $\operatorname{codim}_V W=1$.
	\item[\textup{(b)}] There is a hyperplane $W\subset V$ with $\im A\subseteq W$
	such that $A|_W$ is skew-adjoint for some Euclidean inner product on $W$.
\end{enumerate}

\begin{theorem}\label{th-exist-II}
	Let $\G_A=V\rtimes_A\R e$ be an almost abelian Lie algebra.
	\begin{enumerate}
		\item[\textup{(1)}] $\G_A$ admits a Kundt pair of type~\textup{(II-b)} if and
		only if \textup{(a)} holds. For such pairs $V$ is Lorentzian.
		\item[\textup{(2)}] $\G_A$ admits a Kundt pair of type~\textup{(II-a)} if and
		only if \textup{(b)} holds. In that case, for each of the three signatures
		(Euclidean, degenerate, Lorentzian) there is a Kundt pair of
		type~\textup{(II-a)} such that $\langle\cdot,\cdot\rangle|_V$ has this
		signature.
		\item[\textup{(3)}] Consequently, $\G_A$ admits a Kundt pair of type~\textup{(II)}
		with $V$ Euclidean, or with $V$ degenerate, if and only if \textup{(b)}
		holds; and with $V$ Lorentzian if and only if \textup{(a)} or \textup{(b)}
		holds.
		\item[\textup{(4)}] Condition \textup{(b)} holds if and only if, in a suitable
		basis of $V$,
		\[
		A=\begin{pmatrix}S&b\\0&0\end{pmatrix},\qquad S^T=-S.
		\]
	\end{enumerate}
\end{theorem}

\begin{proof}
	Necessity in (1) and (2). If $p\in W$, Theorem~\ref{th-typeII}(II-b)
	gives $A p\in\R p$, hence the flag $\R p\subset W\subset V$. If $p\notin W$,
	Theorem~\ref{th-typeII}(II-a) gives \textup{(b)}, with the inner product
	induced on $W$.
	
	Sufficiency in (1). Take a basis $(k,w_1,\dots,w_{m-2},t)$ of $V$ adapted to
	the flag, with $L=\R k$ and $W=\spn\{k,w_1,\dots,w_{m-2}\}$, and declare
	\[
	\langle k,t\rangle=1,\quad
	\langle w_i,w_j\rangle=\delta_{ij},\quad
	\langle e,e\rangle=1,
	\]
	all other products being zero. This inner product is Lorentzian, $V$ is
	Lorentzian, and $\h=W\oplus\R e$ is degenerate with $\h^\perp=\R k\subset W$.
	Since $Ak\in\R k$, Theorem~\ref{th-typeII}(II-b) shows that we have a Kundt pair.
	
	Sufficiency in (2) and (4). Assume \textup{(b)}. Take an orthonormal basis
	of $W$ for the given inner product and complete it by some $t\in V\setminus W$.
	The matrix of $A$ is as in (4): $A|_W=S$ is skew-symmetric, and the last row
	vanishes because $\im A\subseteq W$. Conversely, the form in (4) gives
	\textup{(b)} with $W$ the span of the first basis vectors. For $\delta\in\R$,
	define an inner product on $\G_A$ by the Gram matrix
	\[
	G_\delta=
	\begin{pmatrix}
		I_{\dim W}&0&0\\
		0&\delta&1\\
		0&1&0
	\end{pmatrix}
	\quad\text{in the basis }(W,t,e).
	\]
	It is Lorentzian (the block in $(t,e)$ has signature $(1,1)$), and
	$\h=W\oplus\R e$ is degenerate with $\h^\perp=\R e$, $e\notin W$. As
	$\im A\subseteq W$ and $A|_W$ is skew-adjoint, Theorem~\ref{th-typeII}(II-a)
	shows that $(G_\delta,\h)$ is a Kundt pair. The restriction of $G_\delta$
	to $V$ is $I_{\dim W}\oplus(\delta)$: Euclidean for $\delta>0$, degenerate for
	$\delta=0$, Lorentzian for $\delta<0$.
	
	(3) follows from (1), (2) and the fact that $V$ is Lorentzian in
	case~(II-b).
\end{proof}

\begin{remark}
	For $\delta=0$ the same metric makes $V$ a degenerate abelian ideal, so
	it also yields a Kundt pair of type~\textup{(I)} (Theorem~\ref{typeI}).
\end{remark}

\subsection{Consequences for the real Jordan form}

For $\lambda>0$ put
\[
R(\lambda)=\begin{pmatrix}0&-\lambda\\\lambda&0\end{pmatrix},
\qquad J_2(0)=\begin{pmatrix}0&1\\0&0\end{pmatrix}.
\]

\begin{theorem}\label{th-jordan}
	Condition~\textup{(b)} holds if and only if $A$ is similar over $\R$ to
	one of
	\[
	R(\lambda_1)\oplus\cdots\oplus R(\lambda_r)\oplus0_m
	\quad(r\ge1,\ m\ge1),
	\qquad
	R(\lambda_1)\oplus\cdots\oplus R(\lambda_r)\oplus J_2(0)\oplus0_q
	\quad(r,q\ge0),
	\]
	with all $\lambda_j>0$. Equivalently: $0$ is an eigenvalue of $A$, every
	nonzero eigenvalue is purely imaginary and semisimple, and the zero Jordan
	blocks have size at most two, with at most one of size two.
\end{theorem}

\begin{proof}
	By Theorem~\ref{th-exist-II}(4), \textup{(b)} means
	$A=\left(\begin{smallmatrix}S&b\\0&0\end{smallmatrix}\right)$ with
	$S^T=-S$. Decompose $W=\im S\oplus\ker S$ orthogonally and
	$b=b_1+b_0$ accordingly. Replacing $t$ by $t+w$ with $Sw=-b_1$ gives
	$A(t+w)=b_0$. If $b_0=0$, then $A\sim S\oplus0$, which is the first form
	($r\ge1$ because $A\ne0$). If $b_0\ne0$, then $t+w\mapsto b_0\mapsto0$ is a
	nilpotent Jordan chain of length two and the rest of $\ker S$ gives $0_q$;
	this is the second form. In both cases $S$ is orthogonally conjugate to
	$\bigoplus R(\lambda_j)\oplus0$.
	
	Conversely, for the first form let $W$ be the sum of the rotation planes
	and of all but one zero direction; for the second form let $W$ be the
	sum of the rotation planes, of $\ker J_2(0)$ and of $\R^q$. In both cases
	$W$ is a hyperplane containing $\im A$ and $A|_W$ is skew-symmetric for the
	standard inner product. The spectral reformulation is a restatement of
	the two normal forms.
\end{proof}

Let $d_{\R}(A)$ be the sum of the algebraic multiplicities of the real
eigenvalues of $A$.

\begin{theorem}\label{th-lorentz-jordan}
	$\G_A$ admits a Kundt pair for which $V$ is Lorentzian if and only if
	\[
	d_{\R}(A)\ge2
	\qquad\text{or}\qquad
	A\sim[0]\oplus R(\lambda_1)\oplus\cdots\oplus R(\lambda_r),\quad
	r\ge1,\ \lambda_j>0.
	\]
\end{theorem}

\begin{proof}
	A Kundt pair with $V$ Lorentzian is of type~\textup{(II)}, since type~\textup{(I)}
	forces $V$ degenerate. By Theorem~\ref{th-exist-II}(3) it exists if and only
	if \textup{(a)} or \textup{(b)} holds.
	
	\textup{(a)} $\Leftrightarrow d_{\R}(A)\ge2$. Indeed, an invariant flag
	requires a real eigenvector $k$ and a real eigenvalue of the induced map
	on $V/\R k$, and these are two real roots of the characteristic polynomial.
	Conversely, if $d_{\R}(A)\ge2$, choose a real eigenvector $k$; the
	characteristic polynomial of the induced map on $V/\R k$ still has a real
	root, which gives an invariant hyperplane containing $k$.
	
	Assume $d_{\R}(A)<2$ and \textup{(b)}. By Theorem~\ref{th-exist-II}(4) and the
	semisimplicity of $S$, $d_{\R}(A)=1+\dim\ker S$, so $S$ is invertible and
	$b\in\im S$. Taking $t+w$ with $Sw=-b$ as last basis vector gives
	$A\sim S\oplus[0]$, i.e.\ the second alternative with $r\ge1$ (as $A\neq0$).
	Conversely, this normal form satisfies \textup{(b)}.
\end{proof}

\begin{corollary}\label{cor-nilp}
	Let $\G_A$ be nilpotent. Then $\G_A$ admits a Kundt pair for which $V$ is
	Lorentzian. It admits a Kundt pair of type~\textup{(II-a)} (with $V$ of any
	signature) if and only if $A\sim J_2(0)\oplus0_q$, i.e.\
	$\G_A\simeq\mathfrak h_3\oplus\R^q$.
\end{corollary}

\begin{proof}
	For nilpotent $A\neq0$ we have $\dim V\ge2$ and $d_{\R}(A)=\dim V$, so
	Theorem~\ref{th-lorentz-jordan} applies. There are no rotation blocks, so
	Theorem~\ref{th-jordan} leaves only $J_2(0)\oplus0_q$.
\end{proof}

\begin{corollary}\label{cor-indec}
	Let $\G_A$ be indecomposable. Then it admits a Kundt pair of type~\textup{(II-a)}
	if and only if
	\[
	A\sim J_2(0)\oplus R(\lambda_1)\oplus\cdots\oplus R(\lambda_r),
	\qquad r\ge0,\ \lambda_j>0.
	\]
	In particular $\dim\G_A=2r+3$ is odd.
\end{corollary}

\begin{proof}
	Recall that $\G_A$ is indecomposable if and only if $\ker A\subseteq\im A$.
	Among the normal forms of Theorem~\ref{th-jordan}, this holds exactly when
	the only zero Jordan block is $J_2(0)$.
\end{proof}

\begin{remark}
	The preceding results characterize the existence of Kundt pairs for each
	type and each signature of $V$. They do not classify the pairs up to
	isometric Lie algebra automorphisms, which requires a separate analysis
	of the admissible metrics and of their automorphism orbits.
\end{remark}
	
	\section{Kundt Structures of pp-Wave and Plane-Wave Type}
\label{sec:ppwaves}

In this section we identify those Kundt pairs on almost abelian Lie
algebras whose associated left-invariant Lorentzian metrics are pp-waves.
The characterization holds in arbitrary dimension and refines the algebraic
normal forms obtained in Section~\ref{k-a-a}.
We adopt the curvature convention
\(R(X,Y)Z=\nabla_X\nabla_YZ-\nabla_Y\nabla_XZ-\nabla_{[X,Y]}Z\).
The geometric terminology follows
\cite{M-Z,A,G-R-V}; the algebraic description
of Kundt pairs follows \cite{B-M-Z}.

\subsection{Recurrent null lines, Brinkmann spaces and waves}

\begin{Def}\label{def:recurrent-null-line}
	A null line distribution \(\mathcal L=\R K\) on a Lorentzian
	manifold \((M,g)\) is \emph{parallel} (or \emph{recurrent}) if
	\[
	\nabla_XK=\vartheta(X)K\qquad\text{for all }X\in TM
	\]
	for some one-form \(\vartheta\). A nonvanishing null vector field
	\(K\) is \emph{parallel} if \(\nabla K=0\).
	A Lorentzian manifold admitting a parallel null vector field is
	called a \emph{Brinkmann space}. We call the manifold
	\emph{locally Brinkmann} if such a field exists locally.
\end{Def}

The existence of a parallel null line need not imply that its chosen
generator is parallel. Indeed, \(\nabla K=\vartheta\otimes K\)
can be changed to \(\nabla(fK)=0\) locally precisely when
\(\md\vartheta=0\). The recurrence one-forms arising below are
closed and are exact on the simply connected almost abelian group.

\begin{Def}\label{def:pp-plane-waves}
	A \emph{pp-wave} is a Lorentzian manifold possessing a locally
	parallel null vector field \(K\) such that
	\begin{equation}\label{eq:pp-curvature}
		R(X,Y)=0\qquad(X,Y\in K^\perp).
	\end{equation}
	A pp-wave is called a \emph{plane wave} if, in addition,
	\begin{equation}\label{eq:plane-wave-curvature}
		(\nabla_XR)=0\qquad(X\in K^\perp).
	\end{equation}
	A pp-wave is \emph{non-plane} if condition
	\eqref{eq:plane-wave-curvature} fails.
\end{Def}

The curvature formulation is equivalent to the usual Brinkmann
coordinate description; see \cite{M-Z,G-R-V}.
All curvature conditions are first checked relative to the
distinguished null line \(\h^\perp\). The necessity of the conditions
in Theorem~\ref{th:kundt-pp-classification} is intrinsic to the metric:
the scalar curvature vanishes for any pp-wave, and the Ricci tensor
of a pp-wave has rank at most one. We use these obstructions in the
proof, without assuming a priori that another pp-wave direction
coincides with \(\h^\perp\).

\subsection{Algebraic characterization of Kundt pp-waves}

\begin{theorem}\label{th:kundt-pp-classification}
	Let \((\G_A,g,\h)\) be a Kundt pair on an almost
	abelian Lie algebra and let \(G_A\) be the associated simply
	connected Lie group. The metric \(g\) is a pp-wave if and only if
	one of the following conditions holds, according to the type of
	\(\h\):
	\begin{enumerate}
		\item[\textup{(I)}] \(\h=V\). Choose a Witt decomposition
		\[
		V=\R p\oplus E,\qquad
		g(e,p)=1,\quad g|_E=\langle\cdot,\cdot\rangle_E,
		\]
		with \(E\) Euclidean and \(e,p\) null. Write
		\begin{equation}\label{eq:pp-typeI-A}
			A=\begin{pmatrix}a&\alpha\\ b&B\end{pmatrix},
			\qquad b\in E.
		\end{equation}
		Then \(g\) is a pp-wave if and only if \(b=0\), or
		equivalently \(Ap\in\R p\). In that case \(g\) is a
		\emph{plane wave}.
		
		\item[\textup{(II-a)}] \(\h=W\oplus\R e\) and its null
		generator is transverse to \(W\). In an adapted basis \((W,t,e)\),
		\begin{equation}\label{eq:pp-typeIIa}
			A=\begin{pmatrix}S&b\\0&0\end{pmatrix},\quad S^T=-S,
			\qquad
			[g]=\begin{pmatrix}I&0&0\\0&\delta&1\\0&1&0\end{pmatrix},
			\quad\delta\in\R.
		\end{equation}
		Then \(g\) is a pp-wave if and only if \(b=0\); equivalently,
		\(g\) is \emph{flat}.
		
		\item[\textup{(II-b)}] \(\h=W\oplus\R e\) and its null
		generator \(p\) belongs to \(W\). After choosing a generator of the complement that is orthogonal
		to \(V\) and rescaling it to unit spacelike length, choose a Witt basis
		\(V=\R p\oplus E\oplus\R q\) with
		\[
		W=\R p\oplus E,\quad g(p,q)=g(e,e)=1,\quad
		E\perp\operatorname{span}\{p,q,e\},
		\]
		where \(E\) is Euclidean. In this basis
		\begin{equation}\label{eq:pp-typeIIb-A}
			A=\begin{pmatrix}
				a&\alpha&\gamma\\
				0&B&u\\
				0&0&c
			\end{pmatrix}.
		\end{equation}
		Let
		\[
		T=\tfrac12(B+B^*),\qquad
		L=\tfrac12(B-B^*),\qquad
		v=\tfrac12(u+\alpha^\sharp).
		\]
		Then \(g\) is a pp-wave if and only if
		\begin{equation}\label{eq:pp-typeIIb-conditions}
			c=-a,\qquad T=0,\qquad (L+a\Id_E)v=0.
		\end{equation}
	\end{enumerate}
\end{theorem}

\begin{proof}
	We treat the three possibilities separately, deriving the
	connection and curvature formulas needed in each case.
	
	\medskip
	\noindent\textbf{Type (I): the degenerate abelian ideal.}
	Write \(S=(B+B^*)/2\) and \(K=(B-B^*)/2\).
	The brackets are \([e,p]=ap+b\) and
	\([e,x]=\alpha(x)p+Bx\) for \(x\in E\).
	The Koszul formula gives, for $x,y\in E$,
	\begin{align*}
		\nabla_p p&=0,
		&\nabla_p x=\nabla_xp&=\tfrac12\langle b,x\rangle p,\\
		\nabla_p e&=-\tfrac12 b,
		&\nabla_e p&=ap+\tfrac12 b,\\
		\nabla_x y&=\langle Sx,y\rangle p,
		&\nabla_x e&=-Sx-\tfrac12\langle b,x\rangle e,\\
		\nabla_e x&=Kx+\alpha(x)p-\tfrac12\langle b,x\rangle e,
		&\nabla_e e&=-\alpha^\sharp-ae.
	\end{align*}
	Here $\alpha^\sharp$ is the metric dual of $\alpha$ in $E$.
	Tracing the curvature yields the following possibly nonzero
	Ricci components, in addition to symmetry:
	\begin{align*}
		\operatorname{Ric}(p,e)&=\tfrac12\|b\|^2,
		&\operatorname{Ric}(x,y)&=-\tfrac12\langle b,x\rangle\langle b,y\rangle,\\
		\operatorname{Ric}(e,x)&=\tfrac12\langle(B^*+(\operatorname{tr}B)\Id_E)b,x\rangle,
		&\operatorname{Ric}(e,e)&=a\operatorname{tr}S-\operatorname{tr}(S^2)-\alpha(b).
	\end{align*}
	Since the inverse Witt metric pairs \(p\) with \(e\),
	\(\operatorname{Scal}=2\operatorname{Ric}(p,e)+
	\operatorname{tr}_E\operatorname{Ric}|_{E\times E}\). Thus
	\begin{equation}\label{eq:pp-scal-typeI}
		\operatorname{Scal}=\tfrac12\|b\|^2.
	\end{equation}
	The scalar curvature of a pp-wave vanishes, hence \(b=0\)
	is necessary. If \(b=0\), Proposition~\ref{p-recurrent}
	gives
	\[
	\nabla_Xp=a\,\varepsilon(X)p,
	\qquad \varepsilon(e)=1,\quad \varepsilon(V)=0.
	\]
	Here \(\md\varepsilon=0\) because \([\G_A,\G_A]\subset V\).
	Writing \(s\) for the canonical coordinate of the
	\(\R\)-factor on \(G_A\), the field
	\(P=e^{-as}p^L\) is globally parallel.
	
	If $b=0$, the only nonzero curvature terms, up to symmetries, are
	\begin{align*}
		R(e,x)e&=Dx,
		&R(e,x)y&=-\langle Dx,y\rangle p,
		\\[-1mm]
		D&=S^2+[S,K]-aS.
	\end{align*}
	In particular, $R(V,V)=0$. Moreover, $\nabla_xe=-Sx\in V$,
	$\nabla_xy=\langle Sx,y\rangle p$ and $\nabla_xp=0$ for $x,y\in E$;
	all curvature terms involving $p$ vanish. Substitution in the
	covariant derivative of $R$ gives $\nabla_vR=0$ for every $v\in V$.
	This proves that the metric is a plane wave.
	
	\medskip
	\noindent\textbf{Type (II-a): the transverse null direction.}
	The adapted form \eqref{eq:pp-typeIIa} is exhaustive:
	replace the original complement generator by the normalized
	null generator \(p=e+u\), with \(u\in W\), and denote it again
	by \(e\). Since \(W\) is Euclidean, choose
	\(t\in V\cap W^\perp\) with \(g(e,t)=1\).
	Such a choice is possible because \(g(e,\cdot)|_V\)
	has kernel \(W\). Put \(\delta=g(t,t)\) and take an
	orthonormal basis of \(W\). Then \([e,x]=Sx\) and
	\([e,t]=b\), as required.
	The Koszul formula gives the following products.
	For $x,y\in W$ the Levi-Civita product is determined by
	\begin{align*}
		\nabla_x y&=0,
		&\nabla_x t=\nabla_t x&=\tfrac12\langle b,x\rangle(t-\delta e),\\
		\nabla_x e&=-\tfrac12\langle b,x\rangle e,
		&\nabla_e x&=Sx-\tfrac12\langle b,x\rangle e,\\
		\nabla_t e&=-\tfrac12b,
		&\nabla_e t&=\tfrac12 b,\\
		\nabla_t t&=0,&\nabla_e e&=0.
	\end{align*}
	Tracing the curvature, the components needed for the scalar curvature are
	\[
	\operatorname{Ric}(x,y)
	=-\tfrac12\langle b,x\rangle\langle b,y\rangle,\qquad
	\operatorname{Ric}(t,e)=\tfrac12\|b\|^2,
	\qquad \operatorname{Ric}(e,e)=0.
	\]
	The inverse of the \((t,e)\) metric block is
	\(\left(\begin{smallmatrix}0&1\\1&-\delta\end{smallmatrix}\right)\).
	Consequently,
	\begin{equation}\label{eq:pp-scal-typeIIa}
		\operatorname{Scal}=\tfrac12\|b\|^2.
	\end{equation}
	Consequently a pp-wave forces \(b=0\). In this case
	\(A=S\oplus[0]\), and the only possibly nonzero
	Levi-Civita products in the adapted basis are
	\(\nabla_e x=Sx\) for \(x\in W\). Their curvature is
	identically zero. Thus the metric is flat and in particular
	is a pp-wave.
	
	\medskip
	\noindent\textbf{Type (II-b): the null line within the ideal.}
	The adapted normalization \(e\perp V\), \(g(e,e)=1\)
	is achieved by first replacing \(e\) with \(e+w\) for
	\(w\in W\), and then rescaling \(e\) by a nonzero scalar.
	The first change preserves \(A=\ad_e|_V\), while the rescaling
	replaces \(A\) by a nonzero scalar multiple. Thus the normal
	form is not asserted with a fixed numerical matrix \(A\).
	For $v,w\in V$ define the $g|_V$-adjoint $A^\dagger$ and put
	\[
	H=\tfrac12(A+A^\dagger),\qquad
	J=\tfrac12(A-A^\dagger).
	\]
	The Koszul formula gives
	\begin{equation}\label{eq:app-Lorentz-Koszul}
		\nabla_vw=g(Hv,w)e,\qquad
		\nabla_ve=-Hv,\qquad
		\nabla_ev=Jv,\qquad
		\nabla_ee=0.
	\end{equation}
	Contracting the curvature of these connections yields
	\begin{align}\label{eq:app-Lorentz-Ricci}
		\operatorname{Ric}(v,w)
		&=-g\bigl(\bigl([H,J]+(\operatorname{tr}H)H\bigr)v,w\bigr),\notag\\
		\operatorname{Ric}(e,v)&=0,\qquad
		\operatorname{Ric}(e,e)=-\operatorname{tr}(H^2),
	\end{align}
	so that
	\begin{equation}\label{eq:app-Lorentz-scal}
		\operatorname{Scal}=-(\operatorname{tr}H)^2-\operatorname{tr}(H^2).
	\end{equation}
	For the Witt basis $V=\R p\oplus E\oplus\R q$ and the matrix
	$A$ from \eqref{eq:pp-typeIIb-A}, set
	\[
	h=\frac{a+c}{2},\quad k=\frac{a-c}{2},\quad
	v=\frac{u+\alpha^\sharp}{2},\quad
	z=\frac{u-\alpha^\sharp}{2}.
	\]
	Directly from the Witt metric,
	\[
	H=\begin{pmatrix}h&v^*&\gamma\\0&T&v\\0&0&h\end{pmatrix},
	\qquad
	J=\begin{pmatrix}k&-z^*&0\\0&L&z\\0&0&-k\end{pmatrix}.
	\]
	Here \(\operatorname{tr}H=a+c+\operatorname{tr}T\) and
	\(\operatorname{tr}(H^2)=\tfrac12(a+c)^2+\operatorname{tr}(T^2)\).
	Substituting in \eqref{eq:app-Lorentz-scal} gives
	\begin{equation}\label{eq:pp-scal-typeIIb}
		\operatorname{Scal}
		=-(a+c+\operatorname{tr}T)^2
		-\operatorname{tr}(T^2)-\tfrac12(a+c)^2.
	\end{equation}
	If the metric is a pp-wave, its scalar curvature is zero.
	Since \(E\) is Euclidean, each summand on the right-hand
	side is nonpositive. Therefore \(c=-a\) and \(T=0\).
	Under these restrictions, put
	\(v=(u+\alpha^\sharp)/2\) and
	\(z=(u-\alpha^\sharp)/2\). The potentially nonzero
	Ricci components, apart from symmetry, are
	\begin{equation}\label{eq:pp-typeIIb-ricci}
		\operatorname{Ric}(q,x)=\langle (L+a\Id_E)v,x\rangle,
		\quad
		\operatorname{Ric}(q,q)=2a\gamma-2\langle v,z\rangle,
		\qquad x\in E.
	\end{equation}
	All Ricci components with both arguments in the Euclidean
	space \(E\oplus\R e\) vanish. If \(w:=(L+a\Id_E)v\) is nonzero,
	then, for \(x=w\), the restriction of \(\operatorname{Ric}\)
	to \(\operatorname{span}\{q,x\}\) has matrix
	\[
	\begin{pmatrix}\operatorname{Ric}(q,q)&\|w\|^2\\
		\|w\|^2&0\end{pmatrix},
	\] 
	whose determinant is \(-\|w\|^4\ne0\).
	This contradicts the rank-one bound for the Ricci tensor
	of any pp-wave, independently of its null direction.
	Consequently \((L+a\Id_E)v=0\).
	
	Conversely, suppose all the conditions in
	\eqref{eq:pp-typeIIb-conditions} hold. The nonzero
	Levi-Civita products, with all other products vanishing unless explicitly listed, are
	\begin{align*}
		\nabla_xq=\nabla_qx&=\langle v,x\rangle e,&
		\nabla_xe&=-\langle v,x\rangle p,\\
		\nabla_qq&=\gamma e,&
		\nabla_qe&=-\gamma p-v,\\
		\nabla_ep&=ap,&
		\nabla_ex&=Lx-\langle z,x\rangle p,\\
		\nabla_eq&=z-aq,
	\end{align*}
	for \(x\in E\); all products involving \(p\), except
	\(\nabla_ep\), vanish. In particular, the null line
	\(\R p\) is recurrent, with closed recurrence form
	\(a\varepsilon\), and \(e^{-as}p^L\) is parallel on the
	simply connected group. Moreover, the curvature
	components with both first arguments in
	\(\h=\R p\oplus E\oplus\R e\) are determined by
	\[
	R(e,x)q=\langle (L+a\Id_E)v,x\rangle e,
	\qquad
	R(e,x)e=-\langle (L+a\Id_E)v,x\rangle p.
	\]
	All other components with both first arguments in \(\h\)
	vanish, up to skew-symmetry in the first two arguments.
	Hence \(R(\h,\h)=0\) under \eqref{eq:pp-typeIIb-conditions},
	and the metric is a pp-wave.
	
\end{proof}

\subsection{Distinguishing plane waves from non-plane pp-waves}

\begin{pro}\label{prop:typeIIb-plane}
	In the situation of Theorem~\ref{th:kundt-pp-classification}
	\textup{(II-b)}, assume \eqref{eq:pp-typeIIb-conditions}
	and put
	\[
	z=\tfrac12(u-\alpha^\sharp),\qquad
	\kappa=\|v\|^2+2\langle v,z\rangle-2a\gamma.
	\]
	Then the curvature is determined, up to the usual
	symmetries, by
	\begin{align}\label{eq:pp-typeIIb-curvature}
		R(x,q)y&=-\langle v,x\rangle\langle v,y\rangle p,
		&R(x,q)q&=\langle v,x\rangle v,\notag\\
		R(e,q)e&=\kappa p,
		&R(e,q)q&=-\kappa e,
	\end{align}
	for \(x,y\in E\). Furthermore,
	\begin{equation}\label{eq:pp-plane-class}
		g\text{ is a plane wave}\ \Longleftrightarrow\ a\gamma=0.
	\end{equation}
	In particular, the non-plane pp-waves in this class are
	characterized by \(a\gamma\neq0\).
\end{pro}

\begin{proof}
	From the connection formulas in the preceding proof, the
	curvature with first arguments in \(\h\) vanishes, and
	direct substitution into the definition of \(R\) gives
	\eqref{eq:pp-typeIIb-curvature}; all other components vanish
	up to the usual curvature symmetries.
	Since \(L^*=-L\), taking the inner product of
	\((L+a\Id_E)v=0\) with \(v\) gives \(a\|v\|^2=0\).
	For a left-invariant curvature tensor,
	\begin{align*}
		(\nabla_XR)(Y,Z)U={}&\nabla_X(R(Y,Z)U)
		-R(\nabla_XY,Z)U\\
		&-R(Y,\nabla_XZ)U-R(Y,Z)\nabla_XU.
	\end{align*}
	Substitution of the connection and curvature formulas gives
	\begin{equation}\label{eq:app-derivative}
		(\nabla_eR)(e,q)e=2a\kappa p.
	\end{equation}
	If \(a\ne0\), then \(v=0\) and \(\kappa=-2a\gamma\).
	When \(\gamma\ne0\), the derivative in
	\eqref{eq:app-derivative} equals \(-4a^2\gamma p\ne0\).
	Furthermore, \(\operatorname{Ric}(q,q)=2a\gamma\ne0\)
	and the image of the Ricci endomorphism is \(\R p\).
	Any pp-wave direction must therefore coincide with this
	line, so the metric cannot be a plane wave with respect
	to another parallel null direction.
	When \(\gamma=0\), the curvature vanishes identically.
	
	If \(a=0\), then \(Lv=0\). The connection formulas show
	that \(\nabla_x\) sends \(q\) to \(\R e\), \(e\) to
	\(\R p\), and \(E\) to zero for \(x\in E\), whereas
	\(\nabla_e\) acts on \(E\) through \(L\) modulo \(\R p\)
	and sends \(q\) into \(E\). Since
	\(R(E,E)=R(e,E)=R(p,\G_A)=0\) and \(Lv=0\),
	substitution in the covariant derivative formula yields
	\(\nabla_XR=0\) for every
	\(X\in\R p\oplus E\oplus\R e=\h\).
	Thus the metric is a plane wave in this case, completing
	the proof of \eqref{eq:pp-plane-class}.
\end{proof}
The following corollary provides an explicit normal form for Kundt pairs defining non-plane pp-waves, thereby
allowing their construction in arbitrary dimension.
\begin{corollary}\label{cor 4-5}	A Kundt pair on an almost abelian Lie algebra $\G_A$
	defines a non-plane pp-wave if and only if it is of
	type \textup{(II-b)} and, with respect to an adapted
	Witt basis, the defining endomorphism $A$ takes the form
	\[
	A=
	\begin{pmatrix}
		a & \alpha & \gamma\\
		0 & L & -\alpha^\sharp\\
		0 & 0 & -a
	\end{pmatrix},
	\]
	where $L\in\mathfrak{so}(E)$, $\alpha\in E^*$, and
	$a,\gamma\in\R$ satisfy $a\gamma\neq 0$.
\end{corollary}

\begin{proof}
	By Theorem~\ref{th:kundt-pp-classification},
	every Kundt pp-wave of type \textup{(I)} or
	\textup{(II-a)} is a plane wave. Hence a non-plane
	Kundt pp-wave must necessarily be of type
	\textup{(II-b)}.
	
	Combining Theorem~\ref{th:kundt-pp-classification}
	with Proposition~\ref{prop:typeIIb-plane}, we obtain
	\[
	c=-a,\qquad T=0,\qquad
	(L+a\Id_E)v=0,\qquad a\gamma\neq 0.
	\]
	Since $L$ is skew-adjoint, taking the inner product
	of $(L+a\Id_E)v=0$ with $v$ yields
	\[
	0=\langle Lv,v\rangle+a\|v\|^2
	=a\|v\|^2.
	\]
	As $a\neq 0$ and $E$ is Euclidean, it follows that
	$v=0$. Recalling that
	\[
	v=\frac12(u+\alpha^\sharp),
	\]
	we deduce that $u=-\alpha^\sharp$. Moreover, the
	condition $T=0$ implies that $B=L$, with $L^*=-L$.
	Therefore, $A$ has the asserted form.
	
	Conversely, suppose that $A$ has the above form,
	with $L^*=-L$ and $a\gamma\neq 0$. Then
	\[
	c=-a,\qquad T=0,\qquad v=0.
	\]
	Thus all the conditions of
	Theorem~\ref{th:kundt-pp-classification}
	\textup{(II-b)} are satisfied, and the associated
	Lorentzian metric is a pp-wave.
	Since $a\gamma\neq 0$,
	Proposition~\ref{prop:typeIIb-plane} ensures that
	it is not a plane wave.
\end{proof}
The following corollary shows that Ricci-flatness constitutes an obstruction to the existence of non-plane Kundt pp-waves on almost	abelian Lie groups. 
\begin{corollary}\label{co 4-6}
	Every Ricci-flat Kundt pp-wave on an almost abelian
	Lie group is a plane wave.
\end{corollary}

\begin{proof}
	Suppose that a Kundt pair on an almost abelian
	Lie algebra $\G_A$ defines a Ricci-flat pp-wave.
	
	By Theorem~\ref{th:kundt-pp-classification},
	Kundt pp-waves of types \textup{(I)} and
	\textup{(II-a)} are already plane waves.
	It therefore suffices to consider Kundt pp-waves
	of type \textup{(II-b)}.
	
	Assume, by contradiction, that the associated
	Lorentzian metric is not a plane wave.
	By Corollary~\ref{cor 4-5}, we have
	\[
	a\gamma\neq 0,\qquad v=0.
	\]
	On the other hand, the Ricci curvature formula
	\eqref{eq:pp-typeIIb-ricci} yields
	\[
	\operatorname{Ric}(q,q)
	=2a\gamma-2\langle v,z\rangle
	=2a\gamma\neq 0,
	\]
	which contradicts the Ricci-flatness assumption.
	Consequently, the metric is a plane wave.
\end{proof}
\begin{remark}
	Corollary \ref{co 4-6} is consistent with the results of
	Globke and Leistner \cite{G} on
	Ricci-flat locally homogeneous pp-waves.
	Our approach provides a direct algebraic proof
	in the specific setting of almost abelian Lie groups.
\end{remark}

\begin{corollary}\label{cor:ppwave-existence}
	The Lie algebra \(\G_A\) admits a Kundt pair whose
	associated Lorentzian metric is a pp-wave if and only if
	\(A\) has a real eigenvalue. Whenever this is the case,
	there is such a pair of type~\textup{(I)} which is a
	plane wave (possibly flat).
	Consequently every even-dimensional almost abelian Lie
	algebra admits a left-invariant plane-wave metric.
\end{corollary}

\begin{proof}
	If \(Ap=ap\) for some nonzero real vector \(p\), choose
	an adapted Witt metric with \(V^\perp=\R p\).
	Theorem~\ref{th:kundt-pp-classification}\textup{(I)}
	gives a plane wave. Conversely, any of the three cases in
	that theorem provides an eigenvector of \(A\): in
	\textup{(I)} and \textup{(II-b)} it is the adapted
	vector \(p\), whereas in \textup{(II-a)} the condition
	\(b=0\) gives a zero eigenvalue.
	Finally, if \(\dim\G_A\) is even, \(\dim V\) is odd,
	and every real endomorphism of \(V\) has a real
	eigenvalue.
\end{proof}

\subsection{Examples in dimension four and comparison with the literature}

\begin{example}[A non-flat plane wave of type I]
	In the Witt basis \((p,x,y,e)\), equip \(\G_A\) with
	\(g(p,e)=g(x,x)=g(y,y)=1\) and all remaining scalar
	products zero. Let
	\[
	A=\operatorname{diag}(0,1,2)
	\quad\text{on } V=\operatorname{span}\{p,x,y\}.
	\]
	Then \(Ap=0\), so the Kundt pair \((g,V)\) is a
	plane wave. More precisely, the transverse curvature
	operator \(D=S^2+[S,K]-aS\) obtained in the proof of
	Theorem~\ref{th:kundt-pp-classification}\textup{(I)} is
	\(D=\operatorname{diag}(1,4)\), so the metric is not flat.
	In particular, \(\operatorname{Ric}(e,e)=-5\).
\end{example}

\begin{example}[A plane wave of type II]
	Take \(V=\operatorname{span}\{p,x,q\}\) and
	\([e,q]=x\), all other brackets being zero.
	Endow \(\G_A\) with the Lorentzian metric
	\(g(p,q)=g(x,x)=g(e,e)=1\), all other products
	vanishing. The subalgebra
	\(\h=\operatorname{span}\{p,x,e\}\) gives a Kundt pair
	of type~\textup{(II-b)}. Here \(a=\gamma=0\),
	\(v=x/2\), and \(L=0\), so the metric is a
	non-flat plane wave on \(\mathfrak h_3\oplus\R\).
\end{example}

\begin{example}[A non-plane pp-wave of type II]
	In the same Witt basis and metric, let
	\[
	A=\begin{pmatrix}1&0&1\\0&0&0\\0&0&-1\end{pmatrix}.
	\]
	Then \(a=1\), \(\gamma=1\), \(v=0\), and \(L=0\).
	Consequently \((g,\h)\), with
	\(\h=\operatorname{span}\{p,x,e\}\), is a Kundt
	pp-wave which is not a plane wave, since
	\[
	(\nabla_eR)(e,q)e=-4p\ne0.
	\]
	More explicitly, set \(u_1=-q\), \(u_2=-p\),
	\(u_3=x\) and \(u_4=-e\). Then
	\(g(u_1,u_2)=g(u_3,u_3)=g(u_4,u_4)=1\) and
	\[
	[u_1,u_4]=-u_1+u_2,\qquad [u_2,u_4]=u_2,
	\] 
	with all other brackets zero. This is exactly the
	non-plane pp-wave family of\ \cite[Theorem~1.1(1)]{G-R-V}
	with \(\gamma_1=-1\) and \(\gamma_3=0\). In particular,
	\(\G_A\simeq\mathfrak e(1,1)\oplus\R\).
\end{example}

\begin{remark}
	Sato and Tsuyuki \cite[Proposition~3.1 and Lemma~4.4]{S-T}
	studied Ricci-flat metrics on almost abelian groups with a
	degenerate Abelian ideal and proved that the Ricci-flat metrics
	in this family are plane waves. In the type~\textup{(I)} Witt
	coordinates of the present paper, their Ricci-flatness conditions
	become
	\[ b=0,\qquad a\operatorname{tr}S=\operatorname{tr}(S^2),
	\qquad S=\tfrac12(B+B^*). \]
	By contrast, Theorem~\ref{th:kundt-pp-classification}\textup{(I)}
	characterizes the entire plane-wave subfamily by the single
	condition \(b=0\), without assuming vanishing Ricci curvature.
	In fact, the first example above has
	\(\operatorname{Ric}(e,e)=-5\). Thus the Ricci-flat case
	already treated in \cite{S-T} is a proper subfamily of our
	type~\textup{(I)} characterization.
\end{remark}

\begin{remark}
	Theorem~\ref{th:kundt-pp-classification} characterizes
	pp-waves \emph{among the Kundt pairs in adapted normal
		forms}, not their equivalence classes under the entire
	isometry or Lie-algebra automorphism group. In dimension
	four, \cite{G-R-V} gives a classification of
	left-invariant pp-wave metrics on Lie groups up to
	isomorphic homothety. Our statements concern the
	almost abelian subclass, retain the distinguished Kundt
	subalgebra, and hold in arbitrary dimension.
\end{remark}


\begin{thebibliography}{99}
		
		\bibitem{A}
		Aadne, M. T.,
		\emph{Kundt structures},
		arXiv:1912.02570v2 (2019).
		
		\bibitem{Av}
		Avetisyan, Z.,
		\emph{The structure of almost Abelian Lie algebras},
		Int. J. Math. \textbf{33} (2022), no.~8, 2250057.
		doi:10.1142/S0129167X22500574.
		
		\bibitem{B-B}
		Bajo, I., Benayadi, S. and Lebzioui, H.,
		\emph{Flat left-invariant Lorentzian metrics on almost
			abelian Lie groups},
		Results Math. \textbf{81} (2026), 159.
		doi:10.1007/s00025-026-02715-3.
		
		\bibitem{B-M-Z}
		Boucetta, M., Meliani, A. and Zeghib, A.,
		\emph{Kundt three-dimensional left invariant spacetimes},
		J. Math. Phys. \textbf{63} (2022), 112501.
		doi:10.1063/5.0091202.
		
		\bibitem{G-R-V}
		Garc\'{i}a-R\'{i}o, E., Rodr\'{i}guez-Gigirey, R.
		and V\'{a}zquez-Lorenzo, R.,
		\emph{Four-dimensional pp-wave Lie groups and harmonic
			curvature},
		Math. Nachr. \textbf{299} (2026), no.~7, 1751--1780.
		doi:10.1002/mana.70168.
		
		\bibitem{G}
		Globke, W. and Leistner, T.,
		\emph{Locally homogeneous pp-waves},
		J. Geom. Phys. \textbf{108} (2016), 83--101.
		doi:10.1016/j.geomphys.2016.06.013.
		
		\bibitem{M-Z}
		Mehidi, L. and Zeghib, A.,
		\emph{Introduction to Kundt spaces},
		in \emph{Progress in Lorentzian Geometry},
		W. Barrera, J. Herrera, J. P. Navarrete, M. Navarro,
		O. Palmas and D. A. Solis (eds.),
		Springer Proceedings in Mathematics \& Statistics,
		vol.~\textbf{512}, Springer, Cham, 2025, pp.~257--282.
		doi:10.1007/978-3-031-99212-4\_15.
		
		\bibitem{M}
		Meliani, A.,
		\emph{Three-dimensional Lie algebras of Kundt type},
		Eur. Phys. J. Plus \textbf{141} (2026), 897.
		doi:10.1140/epjp/s13360-026-07949-2.
		
		\bibitem{S-T}
		Sato, Y. and Tsuyuki, T.,
		\emph{Lorentzian homogeneous Ricci-flat metrics on
			almost abelian Lie groups},
		J. Math. Phys. \textbf{67} (2026), 112502.
		doi:10.1063/5.0334453.
		
	\end{thebibliography}
\end{document}